\documentclass{artjlt}

\usepackage{amsmath,amsfonts,amssymb}

\usepackage{mathtools}
\input xy
\xyoption{all}

\usepackage{tikz}
\usetikzlibrary{matrix}
\usetikzlibrary{shapes}
\usetikzlibrary{arrows,decorations.pathmorphing,snakes}
\usetikzlibrary{calc,3d}
\usetikzlibrary{decorations,decorations.pathmorphing}
\usetikzlibrary{through}
\tikzset{ext/.style={circle, draw,inner sep=1pt},int/.style={circle,draw,fill,inner sep=1pt},nil/.style={inner sep=1pt}}
\tikzset{exte/.style={circle, draw,inner sep=3pt},inte/.style={circle,draw,fill,inner sep=3pt}}
\tikzset{diagram/.style={matrix of math nodes, row sep=3em, column sep=2.5em, text height=1.5ex, text depth=0.25ex}}
\tikzset{diagram2/.style={matrix of math nodes, row sep=0.5em, column sep=0.5em, text height=1.5ex, text depth=0.25ex}}

\usepackage{scalerel,stackengine}
\usepackage{enumitem}
\stackMath
\newcommand\reallywidehat[1]{\arraycolsep=0pt\relax%
	\begin{array}{c}
		\stretchto{
			\scaleto{
				\scalerel*[\widthof{\ensuremath{#1}}]{\kern-.5pt\bigwedge\kern-.5pt}
				{\rule[-\textheight/2]{1ex}{\textheight}} 
			}{\textheight} %
		}{0.5ex}\\           
		#1\\                 
		\rule{-1ex}{0ex}
	\end{array}
}

\title{Homotopy equivalence of shifted cotangent bundles}                                     
\author{Ricardo Campos}                 
\lastname{Campos}  



\address{%
Ricardo Campos\\               
IMAG, University of Montpellier,
 Place Eugène Bataillon, 
  34090 Montpellier, France\\            
ricardo.campos@umontpellier.fr               
}

\newcommand{\Lie}{\mathsf{Lie}}
\newcommand{\bigslant}[2]{{\raisebox{.2em}{$#1$}\left/\raisebox{-.2em}{$#2$}\right.}}
\newcommand{\Pois}{\mathsf{Pois}}
\newcommand{\Tree}{\mathrm{Tree}}
\newcommand{\id}{\mathrm{id}}

\DeclareMathOperator{\gr}{gr}
\DeclareMathOperator{\sgn}{sgn}
\DeclareMathOperator{\Hom}{Hom}
\DeclareMathOperator{\MC}{MC}
\DeclareMathOperator{\Sh}{Sh}
\DeclareMathOperator{\im}{Im}

\begin{document}


\maketitle

\begin{abstract}
	Given a bundle of chain complexes, the algebra of functions on its shifted cotangent bundle has a natural structure of a shifted Poisson algebra. We show that if two such bundles are homotopy equivalent, the corresponding Poisson algebras are homotopy equivalent. 

We apply this result to $L_\infty$-algebroids to show that two homotopy equivalent bundles have the same $L_\infty$-algebroid structures and explore some consequences about the theory of shifted Poisson structures.   
\end{abstract}

A Lie algebroid consists of a vector bundle $A$ over a manifold $M$ together with a compatible Lie algebra structure on the space of sections $\Gamma(A)$ of $A$. 
More recently, due to the application of homotopy theoretical tools to theoretical physics \cite{kotov2010generalizing,kotov2015characteristic} and to differential geometry (resolution of singular foliations) \cite{lavau2016lie,laurent2017universal}, as well as the study of derived Poisson structures \cite{calaque2017shifted,pridham2017shifted,pantev2018symplectic}, there has been much interest in a derived version of Lie algebroids.

In  the early 90's, T. Lada and J. Stasheff \cite{lada1993introduction} introduced the notion of $L_\infty$ algebras in the context of mathematical physics as a natural extension of differential graded Lie algebras. In an $L_\infty$ algebra, the Jacobi identity is only satisfied up to higher coherent homotopies given by multilinear brackets. 
The same approach of intertwining $L_\infty$ algebras and manifolds gives rise to the homotopical version of Lie algebroids, the so-called $L_\infty$ algebroids \cite{sati2012twisted,severa2001some}.

It is often convenient to work in the dual setting of differential graded (dg) manifolds which are generalizations of smooth manifolds to higher geometry, in which spaces are locally modeled by chain complexes. We recall that in \cite{voronov2010q} Voronov shows that given a graded vector bundle $E$, $L_\infty$ algebroids over $E$ are in one-to-one correspondence with non-positive dg manifold structures on $E$. Given this correspondence, we call $E$ a \textit{split graded manifold}.

Assume now that $E$ dg vector bundle i.e., $E$ is a sequence of vector bundles $(E_i)_{i \in \mathbb Z}$ endowed with a global differential $d\colon E_i \to E_{i+1}$ squaring to zero. 
One of the goals of this manuscript is to understand the behavior of the space of $L_\infty$ algebroid structures on $E$ when we replace $E$ by a homotopy equivalent split dg manifold $F$ over the same base manifold $M$.  

One of our results states that two homotopy equivalent split dg manifolds have essentially the same $L_\infty$ algebroid structures, which can be seen a a version of the Homotopy Transfer Theorem for Lie algebroids, see \cite[Theorem 2.5]{pym2016shifted}.\\

{\textbf{Theorem \ref{thm:algebroid}}}
Let $E$ and $F$ be homotopy equivalent split dg manifolds concentrated in non-positive degrees. Then, there is a bijection

$$\bigslant{\substack{\left\{\substack{\displaystyle L_\infty  \text{ algebroid}\\ \displaystyle \text{structures on } E} \right\}}}{{\text{\small gauge eq.}}} \stackrel{1:1}{\longleftrightarrow} \bigslant{\substack{\left\{\substack{\displaystyle L_\infty  \text{ algebroid}\\ \displaystyle \text{structures on } F} \right\}}}{{\text{\small gauge eq.}}}$$\\

This  correspondence can be obtained by explicit formulas that are given  by sums of trees in the spirit of the homotopy transfer theorem \cite{loday2012algebraic}.
The setting to prove this result is the shifted cotangent bundle $T^*[1]E$ \cite{roytenberg1999courant}. The commutative algebra of functions of this space extends to a shifted Poisson algebra via Kosmann-Schwarzbach's \textit{big bracket} \cite{kosmann2004derived,kosmann1996poisson,voronov2010q}.

There is an of analog of Voronov's result stating that the space of $L_\infty$ algebroid structures over $E$ can be identified with the set of Maurer--Cartan elements of the algebra of functions on $T^*[1]E$, the shifted cotangent bundle of $E$.

This prompts us to understand how the shifted cotangent bundle behaves under homotopy equivalence. Our main result, in the form of Theorem \ref{thm:hom equiv}, states that if $E$ and $F$  are two homotopy equivalent dg vector bundles, their algebras of functions are homotopy equivalent as Poisson algebras.\\

\textbf{Theorem} 
Let $E$ and $F$  be two homotopy equivalent split dg manifolds. Then, there exist $C^\infty(M)$-linear $\infty$-quasi-isomorphisms $\mathcal O_{T^*[1]E} \rightsquigarrow \mathcal O_{T^*[1]F}$ and $\mathcal O_{T^*[1]F} \rightsquigarrow \mathcal O_{T^*[1]E}$ of shifted Poisson algebras. 
Furthermore, this homotopy equivalence of shifted Poisson algebras respects a natural notion of weight.\\

When $E$ is concentrated in degree $0$, $L_\infty$ algebroids are precisely Lie algebroids, and $1$-shifted Poisson structures \cite{calaque2017shifted,pridham2017shifted} are seen to be what is refereed to in the literature as quasi-Lie bialgebroids \cite{bashkirov2016homotopy}. In section \ref{sec:Poisson} we see that under certain conditions our result allows us to conclude that two homotopy equivalent $L_\infty$ algebroids have equivalent spaces of shifted Poisson structures. This matches recent advances by \cite{bandiera2017shifted} and \cite{safronov2017lectures}.\\

\noindent\textbf{Acknowledgments}
I would like to thank Camille Laurent-Gengoux for proposing the problem and applications, as well as many useful discussions. I would also like to thank Damien Calaque for discussions related to Poisson structures and Sylvain Lavau, Pavel Safronov and Joost Nuiten for discussions related to dg geometry. Finally, I wish to thank the anonymous referee for many important comments including finding a mistake with the proof of the main result in the original version of the paper.
I acknowledge support by the Swiss National Science Foundation Early Postdoc.Mobility grant number \texttt{P2EZP2\_174718}.\\

\noindent\textbf{Notation and conventions}
Throughout this manuscript the phrase \textit{differential graded} or \textit{dg} should be implicit everywhere. Concretely, unless otherwise explicit, a vector space $V$ is a dg vector space (i.e. a cochain complex), Lie algebras are differential graded Lie algebras, locally ringed spaces are dg $\mathbb R$-algebras etc. We use cohomological conventions, i.e. all differentials have degree $+1$. In particular this means that taking linear duals negates degrees, that is to say $V^*_i = (V_{-i})^*$.

All vector spaces (such as the ones arising from dg manifolds) considered are assumed to be finite dimensional in every degree and but not necessarily of bounded degree. 

Given two differential graded vector spaces $A$ and $B$, the induced differential on the space $\Hom(A,B)$ is the commutator, denoted by $[d,-]$, satisfying $[d,f] = f\circ d_A + (-1)^{k} d_B\circ f$, for $f \in \Hom(A,B)$ of degree $k$.

The notation $A\rightsquigarrow B$ will be reserved for $\infty$-morphisms of Lie or Poisson algebras, while $A\to B$ will always denote a single map. 

Finally, we consider the ground field to be $\mathbb R$ for concreteness but the reader will notice that all algebraic proofs hold over any field.

\noindent\textbf{Remark about degree shifts}\label{sec:shifts}
Given a vector space $V$, the notation $[k]$ denotes a shift of degree by $k$ units, i.e. $(V[k])_i = V_{k+i}$. 
Throughout the text we will encounter algebraic structures whose operations are not in degree zero. Concretely, the functions on the shifted cotangent bundle form a 2-shifted Lie algebra or a $\Lie\{2\}$ algebra, a Lie algebra whose Lie bracket has degree $-2$. When it is unambiguous, we might omit the shifts for simplicity.

As a precise definition one defines a $\Lie\{k\}$ algebra structure on $V$ to be a Lie algebra structure on $V[-k]$. One should notice that this means that for odd $k$, a $\Lie\{k\}$ algebra has symmetric brackets, but when shifts are even, the defining axioms of (including signs) stay the same.

We remark that one  of the consequences of the degree shifts and the cohomological conventions is that on a $\Lie\{k\}$ algebra, a Maurer--Cartan element has degree $k+1$.

\section{Differential graded manifolds and the shifted cotangent bundle}
In this section, we intend to recall in detail the constructions and results associated to the shifted cotangent bundle of a split dg manifold. We recommend  \cite{fairon2017introduction,Antunes2010,bonavolonta2013category} for a more thorough introduction to the topics of this section.\\

\noindent\textbf{Dg manifolds}

The origins of graded geometry and dg (differential graded) geometry can be traced back to physics, where ($\mathbb Z/2\mathbb Z $ graded) manifolds give for instance a proper treatment of ghosts in BRST deformation. 
Graded (resp. dg) manifolds \cite{kostant} are locally modeled by a graded (resp. dg) vector space $V$ in the sense that a function on such a manifold is locally given by a function on the base manifold and a polynomial function on $V$.
\begin{Definition}
	A graded manifold is a locally ringed space $\mathcal M = (M, \mathcal O_{\mathcal M} )$, where the base $M$ is a
	smooth manifold and around every point $x\in M$ there is an open set $U\ni x$ such that the structure sheaf can be expressed as $\mathcal O_{\mathcal M} (U)= C^\infty(U) \otimes S(V^*)$  for some some graded vector space $V$.
	
	A dg manifold (also called a Q-manifold) a graded manifold equipped with a degree~$+1$ cohomological vector field $Q$, i.e., a derivation of the algebra of functions such that $Q^2=0$.
\end{Definition}

In the present article we will be mostly interested in a subclass of dg manifolds that originate from vector bundles.

\begin{Example}[Dg vector bundles]
	Given a differential graded vector bundle $E$ over $M$, i.e., a sequence of vector bundles ${(E_i)_{i\in \mathbb Z}}$ with differentials $d$
	
	$$\xymatrix{
		...\ E_{i-1} \ar[r]^d \ar[dr] & E_i \ar[d] \ar[r]^d & \ar[dl]E_{i+1}...\\
		& M
	} $$
	
\noindent such that $d^2=0$, one has a naturally associated dg manifold also denoted by $E$, given by its sheaf of sections $E = \left(M,\mathcal O_E = \Gamma(S(E^*))\right)$.
	
	Notice that $d \colon E \to E$ induces a degree $+1$ map $Q\colon E^* \to E^* \subset S(E^*)$ that extends to a square zero $C^\infty(M)$-linear derivation on $\Gamma(S(E^*))$.
	
	Such dg manifolds are called split dg manifolds.
\end{Example}

In fact, Batchelor's theorem \cite{batchelor1980two} (or rather, it's $\mathbb N$-graded version) states that every non-negatively graded manifold originates from such a construction, even though the vector bundle $E$ is non-canonically determined. \\

\noindent\textbf{Shifted cotangent bundle and the big bracket}

Given a graded vector bundle $E\to M$, one can consider its shifted cotangent bundle $T^*[1]E = (M, \mathcal O_{T^*[1]E})$  (see \cite{roytenberg2002structure,Antunes2010} for the constructions in the ungraded setting)\footnote{A more accurate notation for this object from the graded geometry point of view could be $T^*(E[-1])[2]$. In \cite{Antunes2010,roytenberg1999courant} the notation $T^*\Pi E$ is used.} 
Locally this space has coordinates $$x^i \in M, \xi^a\in E, \underbrace{p_i\in TM, \theta_a \in E^*}_{\text{momentum coordinates}}.$$

In these coordinates, the cohomological degree in the algebra $\mathcal O_{T^*[1]E}$ is given by 
$\deg(x^i)=0$, $\deg (p_i)=2$, $\deg (\xi^a) = d+1$ for $\xi^a \in E_d$ and  $\deg(\theta_a) = -d+1$ for $\theta_a\in (E_d)^*$. 
We will also consider a biweight $w$ on $\mathcal O_{T^*[1]E}$ compatible with the product\footnote{In the sense that the product is additive with respect to the biweights.}, where $w(x^i)=(0,0)$, $w(p_i) = (1,1)$, $w(\xi^a) = (0,1)$ and $w(\theta_a)=(1,0)$.

Notice that there are natural inclusions

$$C^\infty|_M \hookrightarrow \mathcal O_{T^*[1]E}, \text{ and }$$
$$
\xymatrix@R-2pc{
	\Gamma(E[-1]) \ar@{^{(}->}[dr] & &\\
 	&  \Gamma(S(E[-1]\oplus E^*[-1])) \ar@{^{(}->}[r]& \mathcal O_{T^*[1]E}.\\
	\Gamma(E^*[-1])\ar@{^{(}->}[ur]	
}$$

\begin{Remark}\label{rem:connection}
	Let us choose connections $\nabla_i$ on $E_i$ for all $i$, and let us consider the corresponding dual connections $\nabla^*_i$ on $E_i^*$. This defines a (non-canonical) inclusion $\Gamma(TM) \hookrightarrow \mathcal O_{T^*[1]E}$. 
	
	With this choice one has an isomorphism of algebras $$\mathcal O_{T^*[1]E} \cong_\nabla S(TM[-2] \oplus E^*[-1] \oplus E[-1]).$$	
\end{Remark}
Besides the commutative product, the space $\mathcal O_{T^*[1]E}$ has a natural Lie bracket  $\{-,-\}$, the so-called \textit{big bracket} \cite{kosmann2004derived,roytenberg2002structure} extending the natural pairing of $E^*$ and $E$.

More concretely, the bracket has degree $-2$, biweight $(-1,-1)$ and it satisfies the following identities on generators

\begin{align*}
\{X,f\} = X\cdot f,& \text{ for } X\in \Gamma(TM), f \in C^\infty|_M,\\
\{\epsilon,e \} = \langle\epsilon,e\rangle,& \text{ for } e \in \Gamma(E), \epsilon\in \Gamma(E^*), \\
\end{align*}

Even though the bracket is intrinsically defined, with the choice of a connection $\nabla$ as in Remark \ref{rem:connection} we also have $\{X,e\} = \nabla_X(e)$ and $\{X,\epsilon \}  = \nabla_X(\epsilon)$.

The bracket is extended to  the full algebra $\mathcal O_{T^*[1]E}$ by the Leibniz rule with respect to the product of functions, making $\mathcal O_{T^*[1]E}$ a  shifted version of a Poisson algebra, also called a $\Pois_3$ or $\mathsf e_3$ algebra in the literature.

\begin{Remark}\label{rem:bidifference}
Since the differential has weight zero and the bracket has weight $(-1,-1)$, the (shifted) Poisson algebra $\mathcal O_{T^*[1]E}$ can be decomposed into a direct sum of (shifted) Lie algebras
	
	$$\mathcal O_{T^*[1]E}= \bigoplus_{k\geq 0} W_{k},$$
	where the Lie algebra $W_k = \bigoplus_{n\geq 0} W_{(n,n+k)}$ is spanned by all the elements whose biweights components have a common difference, i.e, elements of biweight $(0,k), (1,k+1), (2,k+2)$ and so on.	
\end{Remark}

Suppose now that $E$ was a dg vector bundle with differential $d_E$. It is easy to see that these constructions are compatible with the differential and that in this case $\mathcal O_{T^*[1]E}$ is a dg Poisson algebra. 
\begin{Remark}
	
Another way to see this is that $ d_E$ a Maurer--Cartan element of $\mathcal O_{T^*[1]E}$ (seen as a non-differential Poisson algebra), i.e. $\{d_E , d_E\}=0$. Indeed, it follows from $d_E^2=0$ that $\{\{d_E,d_E\}, x\} =0$ for every $x$ element of $E$ or $E^*$. Therefore, $\{d_E,d_E\}$ is central in $S(E\oplus  E^*)$ but the center of this Lie algebra is $\mathbb R$ and therefore $\{d_E,d_E\}=0$.

By twisting the (Lie part of the) Poisson algebra $\mathcal O_{T^*[1]E}$  by this Maurer--Cartan element, we recover a dg Poisson algebra structure on $\mathcal O_{T^*[1]E}^{d_E}$ that we will denote by $\mathcal O_{T^*[1]E}$ only. 

\end{Remark}
\noindent\textbf{(Infinity) Algebroids}
The constructions from the previous section allow us to encode neatly some classical notions. For example, a Lie algebroid structure over $M$ i.e., a Lie algebra bundle $E$ concentrated in degree zero, with a compatible anchor map $\rho \colon E\to TM$, can be expressed as a solution of the Maurer--Cartan equation on $T^*[1]E$:

\begin{Proposition}[\cite{vaintrob1997lie,roytenberg1999courant}]\label{prop:Algebroid is MC}
	Let $M$ be a manifold and $E\to M$ a vector bundle concentrated in degree zero. 
	A Lie algebroid structure on $E$ is equivalent to an element $\mu \in \mathcal O_{T^*[1]E} (M)$ of biweight $(1,2)$ such that $\{\mu,\mu\}=0$.   
	
	The correspondence is given by $\rho(X)\cdot f = \{\{X,\mu\},f\}$ and $[X,Y] = \{\{X,\mu\},Y\}$, for $X,Y \in \Gamma(TM,M)$ and $f\in C^\infty(M)$. 
\end{Proposition}

The same way the \textit{homotopically correct} version of a Lie algebra is an $L_\infty$ algebra, the notion of a Lie algebroid over a manifold $M$ can be homotopically relaxed leading to the concept of an $L_\infty$ algebroid. In what follows we will suppose that all objects are non-positively graded.

\begin{Definition}
	Let $M$ be a smooth manifold and let $(E=(E_i)_{i\leq 0},d)$ be a dg vector bundle over $M$ concentrated in non-positive degree. An $L_\infty$ algebroid structure on $E$ is:
	
	\begin{enumerate}
		\item A dg bundle map $\rho \colon E \to TM$ called \textit{the anchor} and
		\item A sequence of antisymmetric brackets $l_k = [\dots]_k \colon \Gamma(E^{\otimes k}) \to \Gamma(E)$ of degree $2-k$, for $k\geq 2$.
	\end{enumerate}
	
	such that
	
	\begin{enumerate}
		\item All brackets are $C^\infty(M)$ linear except the binary bracket if one of the entries is in degree $0$. If that is the case, then it behaves as a vector field in the sense that if $X\in \Gamma(E_{0})$ and $e\in \Gamma (E)$, 
		$$[X,fe]_2 = f[X,e]_2 +(\rho(X)\cdot f) e.$$
		\item The anchor intertwines $l_2$ and the bracket of vector fields  $$[\rho(x),\rho(y)] = \rho([x,y]), \forall x,y\in \Gamma(E_0).$$
		\item These brackets satisfy the structural axioms of an $L_\infty$ algebra \eqref{eq:Linfty structure}.
	\end{enumerate}
\end{Definition} 

\begin{Remark}
	Some authors such as \cite{getzler2010higher} consider all brackets to be symmetric and of degree $1$ (from an operadic perspective one would call these $L_\infty\{-1\}$ algebroids) while we follow conventions such as the ones of \cite{bonavolonta2013category}. These are equivalent up to a degree shift of $E$.
\end{Remark}

Analogous to Proposition \ref{prop:Algebroid is MC} one can show that $L_\infty$  algebroids are also given as solutions of the Maurer--Cartan equation.

\begin{Proposition}[Folklore]\label{prop:MC=infinity algebroid}
	Let $E\to M$ be a split dg manifold concentrated in non-positive degrees, finite dimensional in every degree. The set of $L_\infty$ algebroid structures over $E$ is in biunivocal correspondence with the space of solutions of the Maurer--Cartan equation in $\mathcal O_{T^*[1]E}$ of biweight $(*,1)$ such that the term in $E^*\otimes E = \Hom(E,E)$ is the differential $d\colon E\to E$.
	
\end{Proposition}

\begin{proof}[Sketch of proof]
	Due to the assumption of finite dimension, a map of bundles $E\to TM$ is equivalent to a section of $E^*\otimes TM$ and the data of the brackets corresponds to a section of $S(E^*)\otimes E$. The degree conditions imply that these correspond to elements of degree $3$ in $\mathcal O_{T^*[1]E}$. 
	
	The bracket condition its easy to verify: The Maurer--Cartan equation can be split by left weight. On left weight $2$ the terms with the differential do not exist due to our degree restraints on $E$. On higher weight we find the $L_\infty$ structure equations and so the Maurer--Cartan equation gives us the same compatibility with the anchor as in the Lie algebroid case.
\end{proof}

\begin{Remark}
	Some authors suppose that $E$ is a graded manifold from the start and the $L_\infty$ algebroid structure includes the datum of the differential $d$ as a unary bracket $l_1$ (see \cite[Definition 1.1.6]{lavau2016lie} for instance). The natural analog of the previous proposition holds, with the differential is recovered from the $E^*\otimes E  $ component. 
	Recall that  the differential $ d_E$ is itself a Maurer--Cartan element of $\mathcal O_{T^*[1]E}$, the two results are related from the general fact that if $\mathfrak g$ is a Lie algebra and $\mu \in \MC(\mathfrak g)$, then $\nu \in \MC(\mathfrak g^{\mu}) \Leftrightarrow \nu+\mu \in \MC(\mathfrak g)$. 
\end{Remark}

\section{Proof of the main result}

The natural notion of homotopy equivalences on cochain complexes generalize naturally to the setting of dg vector bundles.

\begin{Definition}

Two dg vector bundles $E$ and $F$ are said to be homotopy equivalent if there exist bundle maps $f\colon E \to F$ and $g\colon F\to E$ and \textit{homotopies} $H_E\colon E^\bullet \to E^{\bullet+1}$ and $H_F\colon F^\bullet \to F^{\bullet+1}$
	such that $id_E -  g\circ f = H_Ed_E + d_EH_E$ and $id_F -  f\circ g = H_Fd_F + d_FH_F$
	\begin{equation}\label{eq:hom equivalence}	
	\xymatrix{
		E \ar[dr]\ar@(ul,dl)[]_{H_E} \ar@/^0.7pc/[rr]^f
		&& F\ar[dl] \ar@/^0.7pc/[ll]^{g}\ar@(ur,dr)[]^{H_F} \\
		& M & 	
	}
	\end{equation}
\end{Definition}

One can also consider the weaker notion of a \textit{quasi-isomorphism} of dg vector bundles, i.e., a dg vector bundle map $f\colon E \to F$ that induces a quasi-isomorphism on sections. Bear in mind that in general the homology of a dg vector bundle is not a graded vector bundle as it can shift dimensions.

Our main result states that if we take two homotopy equivalent dg vector bundles and consider their shifted cotangent bundles, the respective algebras of functions are homotopy equivalent as Poisson algebras.

\begin{Theorem}\label{thm:hom equiv}
	Let $E$ and $F$ be two dg vector bundles over $M$ that are homotopy equivalent as in the previous definition.

	Then, there exists a $C^\infty(M)$ linear $L_\infty\{2\}$ quasi-isomorphism $\mathcal U\colon \mathcal O_{T^*[1]E} \rightsquigarrow \mathcal O_{T^*[1]F}$.

	Furthermore, this map:
	\begin{itemize}[topsep=0pt,noitemsep]
		\item is compatible with the symmetric algebra product,
		\item is compatible with the biweight in the sense that it preserves each component $W_k$ from Remark \ref{rem:bidifference} (for all $n\geq 1$, $\mathcal U_n$ has biweight $(-n+1,-n+1)$),
		\item its first component $\mathcal U_1$ is the natural extension of $f\oplus g^* \colon E \oplus E^* \to F \oplus F^* $ to a graded commutative algebra morphism.
	\end{itemize}
\end{Theorem}

To be more precise, by compatibility with the symmetric algebra product we mean that every $(\mathcal U_n)_{n\geq 2}$ acts as a derivation with respect to the map $\mathcal U_1$. In particular, this means that $\mathcal U$ actually defines a weak equivalence of shifted Poisson algebras  (an $\infty-\Pois_3$ algebra quasi-isomorphism). This is the notion of morphism considered in \cite{bandiera2017shifted}.

\subsection{The case $M=*$}\label{sec:point}
In this section we prove the main theorem \ref{thm:hom equiv} over $M= *$ a point, which reduces to a problem in homotopical algebra. As we will see, this is the main part of the proof, as the formulas we will obtain over a point readily extend to a more general base.

In this case, $E$ and $F$ are just two dg vector spaces that are quasi-isomorphic with a prescribed homotopy.


The functions on the shifted cotangent bundle $T^*[1]E$ are given by the symmetric algebra $S(E[-1]\oplus E^*[-1])$.

We define a map $\mathcal U_1\colon S(E[-1]\oplus E^*[-1]) \to S(F[-1]\oplus F^*[-1])$ by extending $f\colon E\to F$ and $g^*\colon E^* \to F^*$ to a map of commutative algebras. 

Recall that given a dg vector space $V$, the space $S(V)$ admits a bialgebra structure given by the canonical coproduct $ \Delta \colon S(V) \to S(V)\otimes S(V)$ by $$\Delta(v_1 \dots v_n) = \displaystyle\sum_{\substack{p\leq n \\ \sigma \in \mathbb S_n}} \pm v_{\sigma^{-1}(1)} \dots v_{\sigma^{-1}(p)}\otimes v_{\sigma^{-1}(p+1)}\dots v_{\sigma^{-1}(n)} .$$

Notice that under this description, the Poisson bracket on $S(E[-1]\oplus E^*[-1])=  \mathcal O_{T^*[1]E}$ has the following nice form

$$ \xymatrix@C-1pc{
	& \ar@{-}[d] && \ar@{-}[d]\\
	&\ar@/_0.6cm/@{-}[ddr]\Delta \ar@{-}[dr]&& \ar@{-}[dl]\Delta \ar@/^0.6cm/@{-}[ddl] \\
	l_{ {T^*[1]E}}\coloneqq \{-,-\} =& & \langle -,- \rangle &&, \\
	& & m \ar@{-}[d]&\\
	& & & 
}$$

\noindent where $m$ stands for the multiplication (in the symmetric algebra) and $\langle -,- \rangle $ denotes the pairing between $E^*$ and $E$ being zero otherwise.
By convention, elements of $E^*$ will be placed on the first entry of  $\langle -,- \rangle $ and elements of $E$ will be placed on the second entry.

We define the operator $\mathcal R_2 \coloneqq \mathcal O_{T^*[1]E} \otimes \mathcal O_{T^*[1]E} \to \mathcal O_{T^*[1]E}$

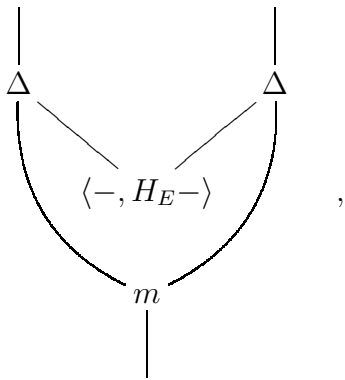
$$ \xymatrix@C-1pc{
	\ar@{-}[d] && \ar@{-}[d]\\
	\ar@/_0.6cm/@{-}[ddr]\Delta \ar@{-}[dr]&& \ar@{-}[dl]\Delta \ar@/^0.6cm/@{-}[ddl] \\
	& \langle -,H_E - \rangle &&, \\
	& m \ar@{-}[d]&\\
	& & 
}$$

Finally, we define $\mathcal U_2\colon \mathcal O_{T^*[1]E} \to\mathcal O_{T^*[1]F}$ to be $\mathcal U_2 \coloneqq  \mathcal U_1 \circ \mathcal R_2$.  

Notice that besides the homotopy, all the operations involved in $\mathcal U_2$ commute with the differentials, from which it follows that 

\begin{equation}\label{eq:U2}
\begin{tikzpicture}
\node at (0,0) {$[d,\mathcal U_2]=$};

\node at (3,0) {
	\xymatrix@C-1pc{
		\ar@{-}[d] && \ar@{-}[d]\\
		\ar@/_0.6cm/@{-}[ddr]\Delta \ar@{-}[dr]&& \ar@{-}[dl]\Delta \ar@/^0.6cm/@{-}[ddl] \\
		& \langle -, - \rangle && \\
		& m \ar@{-}[d]&\\
		& \mathcal U_1 & 
	}
};

\node at (4.6,0) {$-$} ;

\node at (6.5,0) {
	\xymatrix@C-2pc{
		\ar@{-}[d] && \ar@{-}[d]\\
		\ar@/_0.6cm/@{-}[ddr]\Delta \ar@{-}[dr]&& \ar@{-}[dl]\Delta \ar@/^0.6cm/@{-}[ddl] \\
		& \langle -, g\circ f (-) \rangle && \\
		& m \ar@{-}[d]&\\
		& \mathcal U_1 & 
	}
}; 

\node at (10.7,0) {$\eqqcolon \mathcal U_1 \circ l_{ {T^*[1]E}} - \mathcal U_1 \circ \tilde l_{ {T^*[1]E}}$  } ;
\end{tikzpicture}
\end{equation}

\noindent and a similar formula without the terms $\mathcal U_1$ holds if we replace $\mathcal U_2$ by $\mathcal R_2$.
One can easily check on generators that the second term of the equation,  $\mathcal U_1 \circ \tilde l_{ {T^*[1]E}}$ is equal to $\{\mathcal U_1(-),\mathcal U_1(-)\}_{T^*[1]E}$ from where the case $n=2$ from equation \eqref{eq:Linfty morphism} follows.

%

Defining the higher components of the $L_\infty$ morphism requires some set-up. Let $\Tree_n$ be the set of trees with $n$ labeled vertices. To an element $T\in \Tree_n$ one can associate a map $T\colon \mathcal O_{T^*[1]E}^{\otimes n} \to \mathcal O_{T^*[1]E}$ of degree $n-1$. 
The value of $T(x_1,x_2,\dots, x_n)$ is obtained  in the following way:

Let $e_1,\dots, e_{n-1}$ be the set of edges of $T$ and consider a choice of $2n-2$ elements $\alpha_1,\beta_1,\dots,\alpha_{n-1},\beta_{n-1}$ each one of them from either $E$ or $E^*$ such that:

\begin{itemize}
	\item For every $k$, if $e_k$ connects vertices $i$ and $j$, $\alpha_k$ is a factor of $x_i$ and $\beta_k$ is a factor of $x_j$,
	\item There is no repetition of choices.
\end{itemize}

Given such a choice one can consider the product $$\langle \alpha_1,H(\beta_1)\rangle \dots \langle \alpha_{n-1},H(\beta_{n-1})\rangle\stackon[-8pt]{x_1 x_2\dots  x_n}{\vstretch{1.5}{\hstretch{5.0}{\widehat{\phantom{\;}}}}},$$ where  
$ \stackon[-8pt]{x_1 x_2\dots  x_n}{\vstretch{1.5}{\hstretch{5.0}{\widehat{\phantom{\;}}}}}$ denotes the product of all $x_i$'s but with our choice of $\alpha$'s and $\beta$'s removed\footnote{Keep in mind the convention that $\langle a,b \rangle$ is zero unless one of $a,b$ is in $E$ and the other one in $E^*$.} together with the appropriate Koszul sign corresponding to the elements removed.

Finally, the value of $$T(x_1,x_2,\dots, x_n)= \displaystyle\sum_{\substack{\text{choices of }\\ \alpha_1,\dots,\beta_{n-1}}} \langle \alpha_1,H(\beta_1)\rangle \dots \langle \alpha_{n-1},H(\beta_{n-1})\rangle\stackon[-8pt]{x_1 x_2\dots  x_n}{\vstretch{1.5}{\hstretch{5.0}{\widehat{\phantom{\;}}}}}$$ is obtained by summing over all possible choices the products described.

Heuristically, to every edge of a tree we associate an application of the operator $\mathcal R_2$ to its vertices. In particular, $\mathcal R_2= \begin{tikzpicture}[baseline=-.65ex,scale=.6]
\node[ext] (a) at (0,0) {};
\node[ext] (b) at (1,0) {};
\node[below] at (a) {$\scriptstyle 1$};
\node[below] at (b) {$\scriptstyle 2$};
\draw (a) edge (b);
\end{tikzpicture}
\in  \Tree_2 $.

\begin{Remark}\label{rem:sign of graph}
	Notice that some choices regarding the ordering and orientation of edges of $T$ has to be done to compute $T(x_1,x_2,\dots, x_n)$. Since the target $\mathcal O_{T^*[1]E}$ is commutative, all choices lead to the same result up to a sign.
	
	We fix the convention that edges are oriented from the smaller vertex to the bigger vertex and the ordering of edges is done by comparing the smaller label and then the bigger label.
	
	In particular it follows that natural the action of $\mathbb S_n$ permuting the labels of the vertices produces signs.
	
	%
	%
	%
	%
	%
\end{Remark}

For all $n\geq 1$ we define the operators $\mathcal R_n \colon \mathcal O_{T^*[1]E}^{\otimes n} \to \mathcal O_{T^*[1]E}$ of degree $n-1$ as ${\mathcal R_n = \sum_{T \in \Tree_n} T}$. Notice that this definition gives $\mathcal R_1= \id_{\mathcal O_{T^*[1]E}}$.

We also define $\mathcal U_n\coloneqq \mathcal U_1 \circ \mathcal R_n \circ \mathcal O_{T^*[1]E}^{\otimes n} \to \mathcal O_{T^*[1]F}$.

\begin{Proposition}\label{prop: Linfty map over point}
	The maps $\mathcal R_n$ defined above satisfy the following equations, for all $n\geq 2$:	
	
	\begin{equation}
	[d,\mathcal R_n] = \sum_{\sigma \in \Sh^{-1}_{2,n-2}} \sgn(\sigma) \mathcal (\mathcal R_{n-1} \circ_1 l_{T^*[1]E})^\sigma  -  \sum_{\substack{p+q=n\\ \sigma \in \Sh^{-1}_{p,q}}} \sgn(\sigma) (-1)^{p-1} \tilde l_{T^*[1]E}\circ (\mathcal R_p , \mathcal R_q)^\sigma
	\end{equation}
	where $\circ_1$ represents insertion in the first slot and $\tilde l_{T^*[1]E}$ is the twisted bracket defined in equation \eqref{eq:U2}.
\end{Proposition} 

Before proving this proposition notice that by composing the equations above with $\mathcal U_1$ and using the observation that $\mathcal U_1 \circ  \tilde l_{T^*[1]E} =  l_{T^*[1]F}\circ (\mathcal U_1, \mathcal U_1)$ we recover exactly the equations \eqref{eq:Linfty morphism} defining an $L_\infty$ morphism.
\begin{Corollary}
	The maps $\mathcal U_n$ defined above form an $L_\infty$ algebra morphism.	
\end{Corollary}

\begin{proof}[Proof of Proposition \ref{prop: Linfty map over point}] 
	
	Given a tree $T\in \Tree_n$ and $e$ an edge of $T$, we denote by $T^{ e}$ the same tree $T$ but with the edge $e$ replaced by a dashed edge. Similarly, we denote by $T^{\sim e}$ the same edge $e$ replaced by a wavy edge instead.
	
	$$T = \begin{tikzpicture}[baseline=-.65ex,scale=.6]
	\node[ext] (a) at (0,0.5) {};
	\node[ext] (b) at (0,-0.5) {};
	\node[ext] (c) at (1,0) {};
	\node[ext] (d) at (2,0) {};
	\node[ext] (e) at (3,0.5) {};
	\node[left] at (a) {$\scriptstyle 1$};
	\node[left] at (b) {$\scriptstyle 2$};
	\node[below] at (c) {$\scriptstyle 3$};
	\node[below] at (d) {$\scriptstyle 4$};
	\node[below] at (e) {$\scriptstyle 5$};
	\draw (c) edge (a) edge (b) edge (d);
	\draw (d) edge (e);
	\node at (1.5,0.3) {$e$};
	
\end{tikzpicture}
, \quad T^e = 
\begin{tikzpicture}[baseline=-.65ex,scale=.6]
\node[ext] (a) at (0,0.5) {};
\node[ext] (b) at (0,-0.5) {};
\node[ext] (c) at (1,0) {};
\node[ext] (d) at (2,0) {};
\node[ext] (e) at (3,0.5) {};
\node[left] at (a) {$\scriptstyle 1$};
\node[left] at (b) {$\scriptstyle 2$};
\node[below] at (c) {$\scriptstyle 3$};
\node[below] at (d) {$\scriptstyle 4$};
\node[below] at (e) {$\scriptstyle 5$};
\draw (c) edge (a) edge (b);
\draw (d) edge (e);
\draw[dotted] (c) edge (d);
\end{tikzpicture}
, \quad T^{\sim e} = 
\begin{tikzpicture}[baseline=-.65ex,scale=.6]
\node[ext] (a) at (0,0.5) {};
\node[ext] (b) at (0,-0.5) {};
\node[ext] (c) at (1,0) {};
\node[ext] (d) at (2,0) {};
\node[ext] (e) at (3,0.5) {};
\node[left] at (a) {$\scriptstyle 1$};
\node[left] at (b) {$\scriptstyle 2$};
\node[below] at (c) {$\scriptstyle 3$};
\node[below] at (d) {$\scriptstyle 4$};
\node[below] at (e) {$\scriptstyle 5$};
\draw (c) edge (a) edge (b);
\draw (d) edge (e);
\draw[decorate,decoration={snake,amplitude=.7mm,segment length=1mm}] (c) -- (d);
\end{tikzpicture}$$

For $T\in \Tree_n$ we define an action of these modified trees, $T^e,T^{\sim e} \colon \mathcal O_{T^*[1]E}^{\otimes n} \to \mathcal O_{T^*[1]E}$ of degree $n-2$ by the same formula as $T$, except that on the action of the edge corresponding to $e$ connecting vertices $i$ and $j$, with $T^e$ we perform the pairing $\langle \alpha_i,\beta_j\rangle$ and with $T^{\sim e}$ we perform the twisted pairing $\langle \alpha_i,g\circ f (\beta_j)\rangle$.

Notice that since the commutator with the differential $[d,-]$ acts by derivations, the computation of  $[d,\mathcal R_n]$ produces the same kind of composition, except that it replaces one instance of $\langle-,H-\rangle$ by $\langle-,-\rangle- \langle-,g\circ f-\rangle$, just as in equation \eqref{eq:U2}.
In terms of trees, we have that $[d,T] = \sum_{e \text{ edge}} T^e - T^{\sim e}$, so we can also interpret $[d,\mathcal R_n]$ as a sum of all possible trees of $n$ vertices with a dotted edge, minus a sum of all trees with $n$ vertices with a wavy edge.

We claim that the summands corresponding to the terms $T^e$ correspond to the terms $ \sum_{\sigma \in \Sh^{-1}_{2,n-2}} \sgn(\sigma) (\mathcal R_{n-1} \circ_1 l_{T^*[1]E})^\sigma$. 

This follows from the observation that given a tree $\Gamma\in \Tree_{n-1}$, the operation $$T(\{x_1,x_2\},x_3,\dots,x_n)$$ can be expressed as a sum of trees with a dotted edge. 
Concretely, as a quick inspection shows, $\Gamma\circ_1 l_{T^*[1]E}$ is obtained by inserting an graph \begin{tikzpicture}[baseline=-.65ex,scale=.6]
\node[ext] (a) at (0,0) {};
\node[ext] (b) at (1,0) {};
\draw[dotted] (a) edge (b);

\node[below] at (a) {$\scriptstyle 1$};
\node[below] at (b) {$\scriptstyle 2$};
\end{tikzpicture}
on the vertex labeled by $1$ of $\Gamma$ and summing over all possible ways (there exist  $2^{\text{valence of 1}}$) of reconnecting the incident edges, followed by a shift by $1$ of all other labels.

\begin{figure}[h]
$$\begin{tikzpicture}[baseline=1ex,scale=.6]
\node[ext] (a) at (0,0) {};
\node[ext] (b) at (1,0) {};
\node[ext] (c) at (1,1) {};
\draw (a) edge (b) edge (c);

\node[below] at (a) {$\scriptstyle 1$};
\node[below] at (b) {$\scriptstyle 2$};
\node[above] at (c) {$\scriptstyle 3$};
\end{tikzpicture}
\circ_1 \begin{tikzpicture}[baseline=-.65ex,scale=.6]
\node[ext] (a) at (0,0) {};
\node[ext] (b) at (1,0) {};
\draw[dotted] (a) edge (b);

\node[below] at (a) {$\scriptstyle 1$};
\node[below] at (b) {$\scriptstyle 2$};
\end{tikzpicture} = 
\begin{tikzpicture}[baseline=1ex,scale=.6]
\node[ext] (a) at (0,0) {};
\node[ext] (b) at (1,0) {};
\node[ext] (c) at (1,1) {};
\node[ext] (n) at (0,1) {};
\draw (a) edge (b) edge (c);

\draw[dotted] (a) edge (n);

\node[above] at (n) {$\scriptstyle 1$};
\node[below] at (a) {$\scriptstyle 2$};
\node[below] at (b) {$\scriptstyle 3$};
\node[above] at (c) {$\scriptstyle 4$};
\end{tikzpicture}
+
\begin{tikzpicture}[baseline=1ex,scale=.6]
\node[ext] (a) at (0,0) {};
\node[ext] (b) at (1,0) {};
\node[ext] (c) at (1,1) {};
\node[ext] (n) at (0,1) {};
\draw (n) edge (b) edge (c);

\draw[dotted] (a) edge (n);

\node[above] at (n) {$\scriptstyle 1$};
\node[below] at (a) {$\scriptstyle 2$};
\node[below] at (b) {$\scriptstyle 3$};
\node[above] at (c) {$\scriptstyle 4$};
\end{tikzpicture}
+
\begin{tikzpicture}[baseline=1ex,scale=.6]
\node[ext] (a) at (0,0) {};
\node[ext] (b) at (1,0) {};
\node[ext] (c) at (1,1) {};
\node[ext] (n) at (0,1) {};
\draw (a) edge (b) ;
\draw (n) edge (c);

\draw[dotted] (a) edge (n);

\node[above] at (n) {$\scriptstyle 1$};
\node[below] at (a) {$\scriptstyle 2$};
\node[below] at (b) {$\scriptstyle 3$};
\node[above] at (c) {$\scriptstyle 4$};
\end{tikzpicture}
+
\begin{tikzpicture}[baseline=1ex,scale=.6]
\node[ext] (a) at (0,0) {};
\node[ext] (b) at (1,0) {};
\node[ext] (c) at (1,1) {};
\node[ext] (n) at (0,1) {};
\draw (a) edge (c) ;
\draw (n) edge (b);

\draw[dotted] (a) edge (n);

\node[above] at (n) {$\scriptstyle 1$};
\node[below] at (a) {$\scriptstyle 2$};
\node[below] at (b) {$\scriptstyle 3$};
\node[above] at (c) {$\scriptstyle 4$};
\end{tikzpicture}$$
 \caption{An example of an insertion of the graph corresponding to $ l_{T^*[1]E}$.}

\end{figure}
This allows us to conclude that all terms of $(\mathcal R_{n-1} \circ_1 l_{T^*[1]E})^\sigma$ are dotted trees. We just need to show that every dotted tree appears exactly once on the sum over all $(2,n-2)$ unshuffles.

Let us consider an arbitrary dotted tree $T^e$, where $T\in\Tree_n$. Suppose that $e$ connects vertices $i<j$. There is a unique unshuffle $\tau \in \\Sh^{-1}_{2,n-2}$ sending $i$ to $1$ and $j$ to $2$. It is then clear that only $(\mathcal R_{n-1} \circ_1 l_{T^*[1]E})^\tau$ produces trees with a dotted edge connecting vertices $i$ and $j$. Conversely, if we denote by $T/e \in \Tree_{n-1}$ be the graph obtained by the contraction of the edge $e$ one sees that we recover $T^e$ from the insertion of  $l_{T^*[1]E}$ in $T/e$.\footnote{Notice that there is an appearance of a sign factor $\sgn(\tau)$ due to the considerations from Remark \ref{rem:sign of graph}.}\\

To finish the proof it remains to show that $$\sum_{\substack{p+q=n\\ \sigma \in \Sh^{-1}_{p,q}}} \sgn(\sigma) (-1)^{p-1} \tilde l_{T^*[1]E}\circ (\mathcal R_p , \mathcal R_q)^\sigma = \sum_{T\in \Tree_n} \sum_{e \text{ edge of }T} T^{\sim e}.$$

The proof is analogous to the other case. We start by noting that for $T_p\in \Tree_p$ and $T_q\in \Tree_q$,  $\tilde l_{T^*[1]E} \circ (T_i,T_q)$ is obtained summing over all possible ways ($p\times q$) of connecting $T_p$ and $T_q$ with a wavy edge, and shifting the labels of $T_q$ up by $p$ units. It follows that $\tilde l_{T^*[1]E}\circ (\mathcal R_p , \mathcal R_q)^\sigma$ is a sum of elements of the form $T^{\sim e}$, where $T\in \Tree_n$.
To see that every tree appears exactly once, one notices that given a tree with a wavy edge, removing the wavy edge results in a disconnected graph composed of two trees, one in $\Tree_p$ and the other one in $\Tree_{n-p}$ whose labels are uniquely retained by an element of $\Sh^{-1}_{p,q}$.
\end{proof}

\subsection{The global case}

Suppose now that $E$ and $F$ are split dg manifolds over an arbitrary manifold $M$ with maps $f$, $g$ and homotopies $H_E$ and $H_F$ as in equation \eqref{eq:hom equivalence}.

Suppose for the moment being that we can choose connections $\nabla^E$ and $\nabla^F$ that are compatible with $f$ and $g$, i.e. for all $X\in \Gamma(TM), e\in \Gamma(E)$ and $s\in \Gamma(F)$ we have
\begin{equation}\label{eq:connections}
f(\nabla_X^E(e))= \nabla_X^F(f(e)) \text{ and } g(\nabla_X^F(s))= \nabla_X^E(g(s)).
\end{equation}

 Under the identifications induced by $\nabla^E$ and $\nabla^F$,
 $$\mathcal O_{T^*[1]E} \cong S(TM[-2] \oplus E^*[-1] \oplus E[-1]) \text{ and } \mathcal O_{T^*[1]F} \cong S(TM[-2] \oplus F^*[-1] \oplus F[-1]),$$
  the two maps $f$ and $g$ induce a map $\mathcal U_1 \colon \mathcal O_{T^*[1]E} \to \mathcal O_{T^*[1]F}$ of commutative algebras by extending the maps $f\colon E \to F$, $g^* \colon E^* \to F^*$ and $\id \colon TM \to TM$.

Equations \eqref{eq:connections} imply that $\mathcal U_1$ intertwines the Lie brackets on $\mathcal O_{T^*[1]E}$ and $\mathcal O_{T^*[1]F}$ whenever at least one of the elements being bracketed is a vector field. 

We note that the formulas used for $\mathcal U_n, n\geq 2$ in the previous section can be defined over any manifold, since the homotopy $H_E$ is globaly defined on the bundle $E$. It follows that the natural extensions of the maps $\mathcal U_n$ give a well defined $L_\infty$ quasi-isomorphism $ \mathcal O_{T^*[1]E} \rightsquigarrow \mathcal O_{T^*[1]F}$:

\[
\mathcal U_n(x_1,\dots,x_n)\coloneqq\left\{
\begin{array}{ll}
\text{same formula as before if all } x_i \in S(E \oplus E^*) \\
0 \text{ otherwise.}
\end{array}
\right.
\]

Even though is not true in general that we can choose connections $\nabla^E$ and $\nabla^F$ that are compatible with $f$ and $g$, one situation where such choice can be made is if $g\circ f = \id_E$. 

Indeed, the condition $g\circ f = \id_E$  implies that both $f$ and $g$ are maps of constant rank, so their images and kernels are bundles. Identifying $E= \im f$, we can decompose $F = E \oplus \ker g$. 

We can now take an arbitrary connection on $E$, an arbitrary connection on $\ker g$ and define the sum of the two connections as the connection on $F$. This makes the maps $f\colon E \to F$ and $g^* \colon E^* \to F^*$  compatible with the respective connections.

Therefore, the global version of Theorem \ref{thm:hom equiv} follows from the following proposition:

\begin{Proposition}
Given a homotopy equivalence of vector bundles over $M$	$$\xymatrix{
	E \ar@(ul,dl)[]_{H_E} \ar@/^0.7pc/[rr]^f
	&& F \ar@/^0.7pc/[ll]^{g}\ar@(ur,dr)[]^{H_F}	
}$$

there is a dg vector bundle $C$ and homotopy equivalences
$$\xymatrix{
	E  \ar@/^0.7pc/[rr]^{i_E}
	&& C \ar@/^0.7pc/[ll]^{p_E}\ar@(ur,dr)[]^{H_1}	
}
\text{ and }
\xymatrix{
	C \ar@(ul,dl)[]_{H_2} \ar@/^0.7pc/[rr]^{p_F}
	&& F, \ar@/^0.7pc/[ll]^{i_F}	
}
$$

such that $p_E \circ i_E = \id_E$ and $p_F \circ i_F = \id_F$.
\end{Proposition}

\begin{proof}
	We mimic the  standard mapping cylinder construction from homological algebra, see for example \cite{mappingcylinder}. 
	We define $C = E \oplus E[1]\oplus F $ with differential $d(e,e',y) = {(de-e',-de',dy+f(e'))}$.
	
	The second homotopy equivalence depends only on $f$ and is given by the maps $i_F(y) = (0,0,y)$, $p_F(e,e',y) = f(e)+y$ and $H_2(e,e',y) = (0,e,0)$. 
	
	The other homotopy equivalence is given by the maps $i_E(e)= (e,0,0)$, $p_E(e,e',y) = e+H_E(e') + g(y)$ and $H_1(e,e',y) =$ \\
	 
\noindent	 $\begin{pmatrix}
	-g H_F(y+H_Ff(e') -  gfH_E(e')) + H_E 	g(y+H_Ff(e') - fH_E(e') )  + H_E H_E (e')\\
		-g(y+H_Ff(e') - fH_E(e') )  - H_E(e')\\
		 H_F(y+H_Ff(e') + fH_E(e')))	
	\end{pmatrix}.$
\end{proof}

\subsection{Remarks about the hypothesis of homotopy equivalence}

The reader might be surprised that it seems that we almost did not use $F$ (in particular $H_F$) at all in the proofs in section \ref{sec:point}, by reducing the problem to work with $\mathcal R_n$ instead of $\mathcal U_n$.

The reason for this is that one can consider a ``twisted shifted cotangent bundle''  $\widetilde{ \mathcal O_{T^*[1]E}}$ given by the same base space space but with the twisted Lie bracket $\tilde l_{T^*[1]E}$. What we have shown is that $(\mathcal R_n)_{n >1}$ realise an $L_\infty$ isomorphism $\mathcal O_{T^*[1]E}\rightsquigarrow \widetilde{ \mathcal O_{T^*[1]E}}$ extending the identity map.
The result follows from the fact that $\mathcal U_1$ defines a strict Lie quasi-isomorphism $\widetilde{\mathcal O_{T^*[1]E}}\to \mathcal O_{T^*[1]F}$.

However, over manifolds $M$ different from a point, to say that $T^*[1]E$ and $T^*[1]F$ are homotopy equivalent as locally ringed spaces we need the full homotopy data. The reason for this is while $\infty$-quasi-isomorphisms are quasi-invertible over $\mathbb R$, that is not necessarily the case over $C^\infty(M)$.
\section{Applications}

\subsection{Equivalence of $L_\infty$ algebroid structures}

Let $E$ and $F$ be  non-positively graded split dg manifolds that are homotopy equivalent via maps $f$,$g$, $H_E$ and $H_F$ as in the conditions of the main Theorem \ref{thm:hom equiv}. 

Recall from Proposition \ref{prop:MC=infinity algebroid} that $L_\infty$ algebroid structures over $E$ are the same as Maurer--Cartan elements of $\mathcal O_{T^*[1]E}(M)$ of biweight $(*,1)$. 
Since $\mathcal U_n$ has biweight $(-n+1,-n+1)$, it sends $n$ elements of biweight $(*,1)$ to an element of biweight $(*,1)$. It follows that $\mathcal U\colon \mathcal O_{T^*[1]E} \rightsquigarrow \mathcal O_{T^*[1]F}$ maps Maurer--Cartan elements of biweight $(*,1)$ to Maurer--Cartan elements of biweight $(*,1)$. It follows from the Goldman--Millson Theorem \ref{thm:Goldman-Milson} and the main Theorem \ref{thm:hom equiv} that $E$ and $F$ have the same $L_\infty$ algebroid structures.

\begin{Theorem}\label{thm:algebroid}
	Let $E$ and $F$ be split dg manifolds concentrated in non-positive degrees that are homotopy equivalent. Then, there is a set bijection
	%
	%
	
	$$\bigslant{\substack{\left\{\substack{\displaystyle L_\infty  \text{ algebroid}\\ \displaystyle \text{structures on } E} \right\}}}{{\text{\small gauge eq.}}} \stackrel{1:1}{\longleftrightarrow} \bigslant{\substack{\left\{\substack{\displaystyle L_\infty  \text{ algebroid}\\ \displaystyle \text{structures on } F} \right\}}}{{\text{\small gauge eq.}}}$$
	
\end{Theorem}

This result can be compared to the similar result of Pym and Safronov \cite[Theorem 2.5]{pym2016shifted}. Their approach follows the classical proof of the Homotopy Transfer Theorem \cite[Theorem 10.3.9]{loday2012algebraic} while ours is closer to its interpretation in terms of Maurer--Cartan elements and gauge actions as in \cite{dotsenko2016pre}.

\subsection{Isotropy of $L_\infty$ algebroids}

Let us consider $E=(E,d_E,\{l_n^E\}_{n\geq 2})$, an $L_\infty$ algebroid of bounded degree over a manifold $M$ and let us fix a point $m\in M$.

It is well known that on a neighborhood $U\subset M$ of $m$, there exists a dg vector bundle $(F,d_F)$ which is homotopy equivalent to $E|_U$, such that the restriction of the differential $d_F$ to the point $m$ is trivial $d_F|_m = 0$ (see the proof of \cite[Proposition 1.3.5]{laurent2017universal}). 

Using Theorem \ref{thm:algebroid} we can transfer the $L_\infty$ algebroid structure from $E|_U$ to one in $F = (F,d_F,\{l_n^F\}_{n\geq 2})$. 

This structure restricts to an $L_\infty$ algebra on the point $m$, where we have the identification $F_m = H^\bullet(E_m,d_E)$. Notice that when considering the cohomology of an $L_\infty$ algebroid, authors typically consider the anchor $\rho \colon E_0 \to TM$ as part of the cochain complex therefore changing the cohomology in degree zero\footnote{Keep in mind that the anchor is in principle not of constant rank, which means that $\ker \rho \subset F$ is not a vector bundle. This is an important point in the study of singular foliations of \cite{laurent2017universal}.}. In that case, the identification becomes $F_m\supset \ker \rho|_m = H^{\bullet}(E_m,d_E+\rho)$.

Moreover, since the restricted $L_\infty$ algebra structure on $F_m$ has zero differential, $(F_m,l_2^F)$ is a strict (graded) Lie algebra. 
The higher brackets $\{l_n^F\}_{n\geq 3}$ can be seen as  Lie analogs of the Massey products \cite{vallette2014AHO}.

Concretely, the higher brackets $\{l_n^F\}_{n \geq 3} $ correspond to a Maurer--Cartan element in the Chevalley--Eilenberg complex of $(F_m,l_2^F)$, denoted by $\mathrm{CE}(F_m,l_2^F)$. Since the differential is trivial, they actually correspond to an obstruction class living in the Chevalley--Eilenberg cohomology $[l_n^F]\in H^n_{\mathrm{CE}}(H(E_m,d_E),H(E_m,d_E))$.

In fact, since the differential is zero, for every $N\geq 3$, $\{l^F_n\}_{3\leq n\leq N}$ gives an $L_\infty$ structure on $F_m$ and therefore a Maurer--Cartan element in $\mathrm{CE}(F_m,l_2)$. In particular, we obtain the following.

\begin{Proposition}
	Given an $L_\infty$ algebroid $E$ over $M$ and a point $m\in M$, there is a canonically associated class $[l_3] \in H_{\mathrm{CE}}^3(H(E_m,d_E),H(E_m,d_E))$ that vanishes if the $L_\infty$ algebroid structure is homotopically trivial.
\end{Proposition}

The NMRLA (No Minimal Rank Lie Algebroid) class \cite{laurent2017universal} is an example of this class.

\subsection{Shifted Poisson structures}\label{sec:Poisson}

The results that we present here are certainly  connected to the theory of shifted Poisson structures \cite{calaque2017shifted,pridham2017shifted}, see also \cite{bandiera2017shifted,bonechi2018shifted} and \cite{safronov2017lectures}.

Let $E$ be a split dg manifold and $(E,\phi^E)$ an $L_\infty$ algebroid structure on $E$, i.e., $\phi^E$ is a Maurer--Cartan element of $\mathcal O_{T^*[1]E} $, as in the previous section. One can then twist the Poisson algebra $\mathcal O_{T^*[1]E}$ by  $\phi^E$ the Lie algebra $\mathcal O_{T^*[1]E}^{\phi^E} = (\mathcal O_{T^*[1]E}, d^E + \{\phi^E,-\}, \{-,-\})$.
We propose the following definition of a $1$-shifted Poisson structure.

\begin{Definition}
	A $1$-shifted Poisson structure  over the $L_\infty$ algebroid $(E,\phi^E)$ is a Maurer--Cartan element in $\mathcal O_{T^*[1]E}^{\phi^E}$ of biweight $(*,\geq 2)$.
\end{Definition}

This is notion is very close to what in \cite{kravchenko2007strongly,kosmann2005quasi,bashkirov2016homotopy} is referred to as an $L_\infty$ quasi-bialgebra(oid).

From Theorem \ref{thm:hom equiv} and Lemma \ref{lem:twist preserves qi}, we obtain the following result.

\begin{Corollary}
	Let $E$ and $F$ be homotopy equivalent split dg manifolds. Suppose that $\phi^E$ and $\phi^F$ are $L_\infty$ algebroid structures on $E$ and $F$ respectively, such that the map $\mathcal U$ constructed in Theorem \ref{thm:hom equiv} satisfies $\mathcal U(\phi^E) = \phi^F$. 
	Then, the 1-shifted Poisson structures over	$E$ and $F$ are in bijection  up to gauge equivalence.
\end{Corollary}

Recall from \cite{voronov2005higher} ``Voronov's trick'' that out of a decomposition of a Lie algebra $\mathfrak g$ into two Lie sub-algebras $\mathfrak g = \mathfrak h \oplus \mathfrak a$, where $\mathfrak a$ is abelian and a Maurer--Cartan element $\pi \in \MC(\mathfrak h)$ produces an $L_\infty$ structure on $\mathfrak a[1]$ with higher brackets $l_n$ given by the iterated adjoint action $l_n(a_1,\dots,a_n) = \text{pr}_{\mathfrak a}[\dots[\pi,a_1],\dots,a_n]$, where $\text{pr}_{\mathfrak a}\colon \mathfrak g\to \mathfrak a$ denotes the projection. 

Consider an $L_\infty$ algebroid $(E,\phi^E)$. Voronov's trick, applied to $\mathfrak g = \mathcal O_{T^*[1]E}^{\phi^E}$,\ $\mathfrak a = S(\Gamma(E^*[-1]))$,\ $\mathfrak h$ the natural complement (elements of biweight $(*,\geq 1)$) and $\pi$ a $1$-shifted Poisson structure, yields the following proposition: 

\begin{Proposition}
	There exists an $L_\infty$ algebra structure on $S(\Gamma(E^*[-1]))$ whose differential is the one coming from the $L_\infty$ algebroid structure.
\end{Proposition}

Due to the compatibility with the product of functions, we actually obtain that the $L_\infty$ structure extends to a homotopy shifted Poisson structure \cite{calaque2017shifted,pridham2017shifted} that is strict on the product. This is what is called a derived Poisson algebra in \cite{bandiera2017shifted}.

Finally, for $E$ concentrated in degree zero, we recover the classical notion of quasi-Lie bialgebroids (see \cite{roytenberg2002quasi} for a definition of those and, e.g.
\cite{antunes2008poisson},  for the description in terms of big bracket).

\appendix

\section{Recollections about $L_\infty$ algebras and Maurer--Cartan elements}
\newcommand{\Exp}{\mathrm{Exp}}

In this Appendix we recall some of the classical homotopy theory of Lie algebras and their Maurer--Cartan elements that are used in this paper. We assume that the Lie algebras are unshifted, i.e., the bracket has degree zero, but all statements hold for shifted Lie algebras c.f. Section \ref{sec:shifts}.

Recall that an $L_\infty$ algebra structure on the differential graded vector space $(A,d)$ is a family of multilinear antisymmetric maps (the multibrackets) $[-,\dots,-] = l_n \colon A^{\otimes n} \to A$ of degree $|l_n| = 2-n$ for $n\geq 2$ satisfying the higher Jacobi identities: 

\begin{equation}\label{eq:Linfty structure}
\sum_{\substack{p+q=n+1\\ p,q>1}} \sum_{\sigma\in \Sh^{-1}_{q,p-1}}\sgn(\sigma) (-1)^{(p-1)q}(l_p\circ_1 l_q)^{\sigma} = [d,l_n],
\end{equation}
where $\Sh^{-1}_{q,p-1}\subset \mathbb S_{q+p-1}$ denotes the $(q,p-1)$ unshuffles.

Most results in this section can be generalized to $L_\infty$ algebras but since they are not necessary for us they are stated in terms of Lie algebras and $L_\infty$ morphisms for simplicity of formulas.

\begin{Definition}
	An $L_\infty$ morphism $\mathcal U \colon A\rightsquigarrow B$ between two Lie algebras $(A,l_A,d_A)$ and $(B,l_B,d_B)$ is a sequence of maps $\mathcal U_n \colon S^n A \to B, \forall n\geq 1$ of degree $1-n$ such that $\mathcal U_1$ commutes with the differentials, i.e $[d,\mathcal U_1]=0$ and 
	\begin{equation}\label{eq:Linfty morphism}
	[d,\mathcal U_n] = \sum_{\sigma \in \Sh^{-1}_{2,n-2}} \sgn(\sigma) \mathcal (\mathcal U_{n-1} \circ_1 l_A)^\sigma  -  \sum_{\substack{p+q=n\\ \sigma \in \Sh^{-1}_{p,q}}} \sgn(\sigma) (-1)^{p-1} l_B\circ (\mathcal U_p , \mathcal U_q)^\sigma
	\end{equation}
\end{Definition}

\begin{Definition}
	Let $\mathfrak g$ be a differential graded Lie algebra. A \textit{Maurer--Cartan} element is an element $\mu\in \mathfrak g^1$ of degree $1$ that satisfies the equation
	$$d\mu + \frac 1 2[\mu, \mu]=0$$
	
	The set of Maurer--Cartan elements of a Lie algebra $\mathfrak g$ is denoted by $\MC(\mathfrak g)$.
\end{Definition}

\begin{Definition}
	A filtered Lie algebra is a Lie algebra $\mathfrak g$ equipped with a complete descending filtration $\mathcal F^\bullet$ of Lie algebras i.e. $\mathfrak g =\mathcal F^1 \mathfrak g \supset \mathcal F^2 \mathfrak g \supset\mathcal F^2 \mathfrak g \supset \dots$ satisfying $\left[\mathcal F^i \mathfrak g,\mathcal F^j \mathfrak g\right] \subset \mathcal F^{i+j} \mathfrak g$, such that $\mathfrak g$ is complete with respect to this filtration
	$$\mathfrak g = \varprojlim_{k} \mathfrak g / \mathcal F^k \mathfrak g.$$
\end{Definition}

Let $\mathfrak g$ and $\mathfrak h$ be filtered Lie algebras. It is easy to check that given an $L_\infty$ morphism $\mathcal U = (\mathcal U_k)_{k \geq 1} \colon \mathfrak g \rightsquigarrow \mathfrak h$ compatible with the filtrations\footnote{In the sense that $\mathcal U_k(\mathcal F^{i_1}\mathfrak g,\dots,F^{i_k}\mathfrak g) \subset \mathcal F^{i_1+\dots +i_k}\mathfrak h$.} and a Maurer-Cartan element $\mu\in \mathfrak g$, then, the element 
\begin{equation}\label{eq:U(Maurer-Cartan)}
\mathcal U(\mu)\coloneqq\sum_{n=1}^\infty \frac{1}{n!}\mathcal U_n(\mu,\dots,\mu) \in  \varprojlim \mathfrak h / \mathcal F^k \mathfrak h=\mathfrak h
\end{equation}
is a Maurer-Cartan element of $\mathfrak h$.

Given a Maurer-Cartan element $\mu$ of a Lie algebra $\mathfrak g$ one often considers the corresponding twisted Lie algebra $\mathfrak g^\mu$.

\begin{Definition}
	Let $\mathfrak g$ be a differential graded Lie algebra and $\mu\in \MC(\mathfrak g)$. We denote by $\mathfrak g^\mu$ the \textit{twist} of $\mathfrak g$ by $\mu$, which is a differential graded Lie algebra that is equal to $\mathfrak g$ as a graded Lie algebra, with differential given by
	$$d_{\mathfrak g^\mu} = d_\mathfrak g + [\mu,-].$$
\end{Definition}

Twisting is a homotopically stable property. The following result follows from a simple spectral sequence argument.

\begin{Proposition}[\cite{dolgushevthesis}, Proposition 1]\label{lem:twist preserves qi}
	Let $\mathfrak g$ and $\mathfrak h$ be Lie algebras and $U\colon \mathfrak g \to \mathfrak h$ be an $L_\infty$ morphism.
	
	If for all $k$, $\mathcal U_1\colon \mathcal F^k \mathfrak g \to \mathcal F^k \mathfrak h$ is a quasi-isomorphism, then for any $\mu \in \MC(\mathfrak g)$, the induced map $\mathcal U^\mu\colon \mathfrak g^\mu \rightsquigarrow \mathfrak h^{\mathcal U(\mu)}$ is an $L_\infty$ quasi-isomorphism.
\end{Proposition}

Given a Lie algebra $\mathfrak g$ and an commutative algebra $A$, the space $\mathfrak g\otimes A$ inherits a natural Lie algebra structure by declaring the bracket to be $A$-bilinear, i.e., $[X\otimes a, X'\otimes a'] = [X,X'] \otimes aa'$. In the case of the polynomial forms  $A=\Omega_{\text{poly}}([0,1]) = \mathbb K [t,dt]$, we get a natural Lie algebra structure on $\mathfrak g[t,dt]$.

\begin{Definition}	Let $\mathfrak g$ be a Lie algebra. Two Maurer--Cartan elements $\mu_0,\mu_1\in \MC(\mathfrak g)$ are said to be gauge equivalent if there is a Maurer--Cartan element $\mu_t \in \mathfrak g[t,dt]$ interpolating $\mu_0$ and $\mu_1$.
\end{Definition}

This definition amounts to say that $\mu_t$ can be written for all $t\in [0,1]$ as
\[
\mu_t = m_t + h_tdt 
\]
where $m_t$ can be understood as a family of Maurer-Cartan elements in $\mathfrak g$, connected by a family of infinitesimal homotopies (gauge transformations) $h_t\in \mathfrak g^0$.
The Maurer--Cartan equation for $\mu_t$ translates into the two equations
\begin{align*}
dm_t+\frac 1 2 [m_t,m_t]&=0,
&
\dot m_t + dh_t +[h_t,m_t]&=0.
\end{align*}

Remarkably, the Goldman--Millson theorem states that under appropriate conditions one can identify the Maurer--Cartan spaces of quasi-isomorphic Lie algebras.

\begin{Theorem}[Goldman--Millson \cite{dolgushev2015version}]\label{thm:Goldman-Milson}
	Let $U\colon \mathfrak g \to \mathfrak h$ be an $L_\infty$ morphism of filtered Lie algebras. Suppose furthermore that on the associated graded level the map $\gr U\colon \gr \mathfrak g = \bigoplus \mathcal F^\bullet \mathfrak g /\mathcal F^{\bullet+1} \mathfrak g  \to \gr \mathfrak h$ is a quasi-isomorphism. Then, formula \eqref{eq:U(Maurer-Cartan)} induces a bijection of sets $$U\colon \MC(\mathfrak g)/\text{gauge equiv.} \to \MC(\mathfrak h)/\text{gauge equiv.}$$
\end{Theorem}



\end{document}